\documentclass{ifacconf}
\usepackage{graphicx}
\usepackage{mathtools}   
\usepackage{graphicx}    
\usepackage{algcompatible}
\usepackage{algorithm}
\usepackage{enumerate}   

\usepackage{natbib}             
\usepackage{amsfonts}
\usepackage{savesym}
\savesymbol{AND}
\usepackage{xcolor}
\newcommand*{\Scale}[2][4]{\scalebox{#1}{$#2$}}%
\usepackage{bm}
\definecolor{Blue}{RGB}{88, 105, 225}
\definecolor{dred}{rgb}{0.6,0,0}

\def\eqref#1{equation~\ref{#1}}

\def\1{\bm{1}}

\def\rx{{\textnormal{x}}}

\def\vmu{{\bm{\mu}}}

\def\va{{\bm{a}}}

\def\vc{{\bm{c}}}
\def\vd{{\bm{d}}}

\def\vs{{\bm{s}}}

\def\vu{{\bm{u}}}
\def\vv{{\bm{v}}}
\def\vw{{\bm{w}}}
\def\vx{{\bm{x}}}
\def\vy{{\bm{y}}}

\def\mA{{\bm{A}}}
\def\mB{{\bm{B}}}
\def\mC{{\bm{C}}}

\def\mG{{\bm{G}}}

\def\mM{{\bm{M}}}

\def\mP{{\bm{P}}}
\def\mQ{{\bm{Q}}}
\def\mR{{\bm{R}}}
\def\mS{{\bm{S}}}

\def\mSigma{{\bm{\Sigma}}}

\DeclareMathAlphabet{\mathsfit}{\encodingdefault}{\sfdefault}{m}{sl}
\SetMathAlphabet{\mathsfit}{bold}{\encodingdefault}{\sfdefault}{bx}{n}

\def\sF{{\mathbb{F}}}

\def\sH{{\mathbb{H}}}

\def\sN{{\mathbb{N}}}

\def\sP{{\mathbb{P}}}

\def\sS{{\mathbb{S}}}

\newcommand{\E}{\mathbb{E}}

\newcommand{\R}{\mathbb{R}}

\begin{document}
\begin{frontmatter}

\title{
Towards Resilient Tracking in Autonomous Vehicles: A Distributionally Robust Input and State Estimation Approach
}

\thanks[footnoteinfo]{This work is supported by the University of New Mexico under the School of Engineering (SOE) faculty startup funding and under the Research Allocations Committee (RAC) grant \#7gr4Py.}

\author{Kasra Azizi \ \ }
\author{Kumar Anurag \ \ } 
\author{Wenbin Wan} 
\address{University of New Mexico, Albuquerque, NM, USA \\
\{azkasra, kmranrg, wwan\}@unm.edu}

\begin{abstract}
This paper proposes a novel framework for the distributionally robust input and state estimation (DRISE) for autonomous vehicles operating under model uncertainties and measurement outliers. The proposed framework improves the input and state estimation (ISE) approach by integrating distributional robustness, enhancing the estimator's resilience and robustness to adversarial inputs and unmodeled dynamics. Moment-based ambiguity sets capture probabilistic uncertainties in both system dynamics and measurement noise, offering analytical tractability and efficiently handling uncertainties in mean and covariance.
In particular, the proposed framework minimizes the worst-case estimation error, ensuring robustness against deviations from nominal distributions.
The effectiveness of the proposed approach is validated through simulations conducted in the CARLA autonomous driving simulator, demonstrating improved performance in state estimation accuracy and robustness in dynamic and uncertain environments. 
\end{abstract}

\end{frontmatter}
\section{Introduction} \label{sec: intro}
\begin{figure*}
\begin{center}
\includegraphics[width=0.95\textwidth]{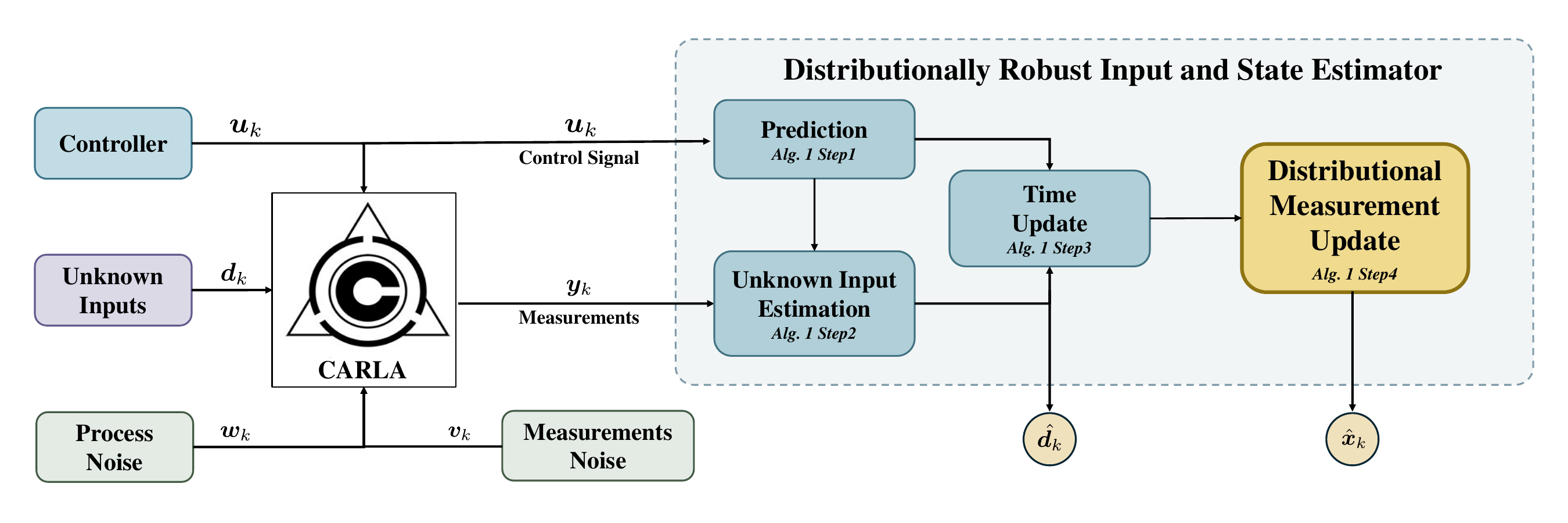}   
\caption{Distributionally Robust Input and State Estimation (DRISE) framework diagram. The estimation process follows the steps outlined in Algorithm~\ref{alg_DRISE}.}
\label{fig: block_diagram}
\end{center}
\end{figure*}

In recent decades, advances in automotive technologies have led to the widespread integration of active safety features in modern vehicles, such as antilock braking systems, electronic stability control, traction control systems, direct yaw moment control, and active suspension systems (\cite{9744542}). These systems significantly improve vehicle handling, stability, and overall safety. A fundamental aspect of their functionality is the accurate and real-time availability of vehicle state information and driver-related variables, which are essential for their operation. Although high-resolution sensors and redundant setups can provide precise data, the cost of such implementations can be prohibitive for mass-production vehicles. Consequently, there is a pressing need for reliable and cost-effective real-time estimation algorithms that can utilize affordable sensors to generate accurate driver-vehicle information, thereby enabling the efficient design of active safety systems for modern vehicles (\cite{safety_systems_new}).

\textit{Literature Review.}
Extensive research has been implemented on accurately estimating vehicle states that are difficult to measure directly with cost-effective sensors, such as longitudinal and lateral velocities (\cite{m_estimation}). Resilient estimation (RE) has emerged to address the challenge of estimating true system states even when measurements are disrupted, ensuring system integrity by reducing the impact of disturbances on estimation quality (\cite{wan2024safeGPS,cps_survey_new}).
 
For deterministic linear systems, secure state estimation is often approached using set-based methods, which provide robustness against adversarial sensor attacks by ensuring the true state remains within a computed estimation set (\cite{set_based_1}).
These methods construct bounded sets based on sensor measurements and system dynamics, mitigating the impact of corrupted measurements.
However, these approaches generally assume noise-free environments and perfectly known dynamics, making it difficult to guarantee resilience under more realistic conditions. For stochastic linear systems, when process and measurement noises are Gaussian and model parameters are exact, the Kalman filter achieves optimal solutions in various senses, such as minimum-variance unbiased estimation and least-squares error minimization (\cite{unified,new_Kitanidis}). However, these traditional methods lack robustness against unexpected disturbances or outliers, which could lead to degraded performance or divergence in state estimates.

Existing resilient estimation methods, such as those for linear (\cite{RE_linear_new}) and nonlinear stochastic systems (\cite{hunmin_NL}), focus on minimizing estimation error variance under adversarial conditions. In particular, they incorporate unknown disturbances as part of the system dynamics, allowing for detection algorithms based on statistical thresholds, such as the widely-used $\mathcal{X}^2$ test. While these methods effectively detect anomalies, they often rely on assumed disturbance models, which can be difficult to specify accurately in real-world scenarios with uncertainties (\cite{wan_attack_resilient}).
Real-world engineering applications often face non-trivial uncertainties in system parameters, such as unknown target maneuvers (\cite{manouvering_new}), sensor faults (\cite{sensor_faults_new}), and noise environment variability and measurement anomalies (\cite{changing_noise_new}). The Kalman filter is particularly sensitive to these deviations, which can deteriorate its performance. Thus, robust estimation methods, such as moment-based distributionally robust estimators (\cite{moment_based,DRE}), have been developed to address parameter uncertainty by minimizing the worst-case estimation error.
Another line of research focuses on making state estimation insensitive to outliers, notably with M-estimation-based Kalman filters (\cite{m_estimation_kf_new}). These methods identify and mitigate the influence of outliers using influence functions that limit their effect on estimation.

 Traditional robust estimation techniques typically rely on detailed structural information about uncertainties, which is often unavailable in practical scenarios. To overcome this issue, \cite{DRE} proposed a distributionally robust estimation (DRE) approach that optimizes over a family of distributions within a specified distance from a nominal distribution. This approach uses ambiguity sets to account for uncertainty, enabling robust state estimation even in the presence of parameter uncertainties and measurement outliers. However, while DRE provides significant advancements in handling distributional uncertainties, it does not explicitly account for unknown inputs in its approach. This limitation restricts its applicability to real-world systems where unmeasured inputs are critical.
 
\textit{Contributions.}  
Motivated by the limitations of distributionally robust estimation (DRE) in handling unknown inputs, this paper introduces the distributionally robust input and state estimation (DRISE) framework as shown in Figure~\ref{fig: block_diagram}. The main contributions of this work are summarized as follows.
\begin{enumerate} [(i)]
    \item The DRISE framework improves the DRE approach by incorporating the estimation of unknown inputs alongside states. This innovation addresses a critical gap in the existing procedures and enables robust performance in systems subjected to unknown inputs.
   
    \item  To the best of our knowledge, we are the first to employ moment-based ambiguity sets into state estimation while taking unknown inputs into account. This ensures robustness against measurement outliers, deviations from nominal noise distributions, and time-varying uncertainties. The proposed framework improves estimation reliability in complex environments.
   
    \item  We demonstrate the practical applicability of the proposed framework in the CARLA autonomous driving simulator. We showcase its enhanced robustness and performance under uncertainties by benchmarking it against the Kalman Filter, ISE, and DRE in a realistic driving environment.
\end{enumerate}

The remainder of this paper is organized as follows. Section~\ref{sec: prelim} introduces the preliminaries, including notations and the system model. Section~\ref{sec: problem} formalizes the problem statement. Section~\ref{sec: DRISE} presents the proposed DRISE framework. Section~\ref{sec: exprements} evaluates the DRISE framework through numerical experiments using CARLA simulation. 
Section~\ref{sec: further} discusses potential future works. Finally, Section~\ref{sec: conclusion} concludes the paper. 

\section{Preliminaries} \label{sec: prelim}

\textit{Notations.} 
Let $k$ denote the time index. We use $\R^n_+$ to denote the set of positive elements in $\R^n$, and $\R^{n \times m}$ for the set of all $n \times m$ real matrices. For a matrix $\mA$, $\mA^\top$, $\mA^{-1}$, $\mA^+$, $\text{diag}(\mA)$, $\text{tr}(\mA)$, and $\text{rank}(\mA)$ represent the transpose, inverse, Moore-Penrose pseudoinverse, diagonal, trace, and rank of $\mA$, respectively.  
A symmetric matrix $\mS \in \sS_d$ satisfies $\mS=\mS^\top$. We write $\mS>0$ ($\mS \geq 0$) to indicate that $\mS$ is positive (semi)definite. Let $\sS_d^+ \text{(resp. } \sS_d^{++}$) be the set of symmetric positive semi-definite (resp. positive definite) matrices. If $\mS \in \sS_d^+$, $\mS^{1/2}$ denotes the square root of $\mS$, i.e., $\mS^{1/2}\mS^{1/2}=\mS$.  
We use $\|\cdot\|$ to denote the Euclidean norm or induced matrix norm. For a vector $\va$, $\va_i$ is the $i$-th element of $\va$.  
Let $\E[\cdot]$ denote the expectation operator with respect to the distribution $\sP$. For a vector $\va$, $\hat{\va}$ and $\tilde{\va} = \va - \hat{\va}$ represent the estimate and estimation error, respectively.  
For distributions, $\mathcal{N}(\vmu, \mSigma)$ is a $d$-dimensional  Gaussian with mean $\vmu$ and covariance $\mSigma$. The joint or conditional distribution of $\{\rx, \dots\}$ is written as $\sP_x(\rx, \dots)$.  
Finally, $\mathbf{Y}_k$ represents the measurement sequence up to time $k$, i.e., $\mathbf{Y}_k = \{\vy_1, \dots, \vy_k\}$.
To minimize notational complexity, an ellipsis within square brackets indicates the repetition of the preceding bracketed expression. For example, $[\xi][\cdots]$ denotes $[\xi][\xi]$, where $[\xi]$ represents a potentially lengthy expression.

\textit{System Model.}
Consider the following linear time-varying (LTV) discrete-time stochastic system 
\begin{subequations} \label{eq1}
\begin{align}
    \vx_{k+1} &= \mA_k \vx_k + \mB_k \vu_k + \mG_k \vd_k + \vw_k  \\
    \vy_k &= \mC_k \vx_k + \vv_k,
\end{align}
\end{subequations}
where $\vx_k \in \R^n$ is the state vector at time $k=1, 2, \hdots$ and $\vu_k \in \R^m$ is a known control input vector. $\vd_k \in \R^p$ is an unknown input vector, and $\vy_k \in \R^l$ is the measurement vector. 
The process noise $\vw_k \in \R^n$ and the measurement noise $\vv_k \in \R^l$ are assumed to be mutually uncorrelated, zero-mean, white random signals with known covariance matrices, $\mQ_k = \E[\vw_k \vw_k^\top] \geq 0$ and $\mR_k = \E[\vv_k \vv_k^\top] > 0$, respectively. 
The matrices $\mA_k, \mB_k, \mG_k$, and $\mC_k$ are known and have finite matrix norms. $\vx_0$ is assumed to be independent of $\vv_k$ and $\vw_k$ for all $k$ and the unbiased estimate $\hat{\vx}_0$ of the initial state $\vx_0$ is available with covariance matrices $\mP^x_0, \mP^{xd}_0$ and $\mP^d_0$. 
In addition, we have the following two assumptions.
\begin{assum}
The system has perfect/strong observability, i.e., the initial condition $\vx_0$ and the unknown input sequence $\{\vd_i\}_{i=0}^{r-1}$ can be uniquely determined from the measured output sequence $\{\vy_i\}_{i=0}^r$ of a sufficient number of observations, i.e., $r \geq r_0$ for some $r_0 \in \sN$ (\cite{payne1973discrete}).
\end{assum}
\begin{assum}    
 We assume that $\text{rank}(\mC_k\mG_{k-1})=p$, where $0 \leq p \leq l$. The interpretation of this assumption is that the impact of the unknown inputs $\vd_{k-1}$ on the system dynamics can be observed by $\vy_k$. This is a typical assumption, as in \cite{CARE_34}.
\end{assum}
 
\section{Problem Statement} \label{sec: problem}
Consider a general estimation problem where the goal is to minimize the estimation error of a state $\vx_k$ based on the observation $\mathbf{Y}_k$. This problem is formulated as    
\begin{equation}\label{eq: DRE_optimization}
    \min_{\phi(\cdot) \in \sH^\prime_{\mathbf{Y}_k}} \E_{\bar{\sP}(\vx_k, \mathbf{Y}_k)} \big[ \vx_k - \phi(\mathbf{Y}_k) \big]  \big[\vx_k - \phi(\mathbf{Y}_k) \big]^\top,
\end{equation}   
where $\sH^\prime_{\mathbf{Y}_k}$ represents the set of all Borel measurable functions of $\{1, \mathbf{Y}_k\}$ with finite second moments, $\phi(\cdot)$ is the estimator, and $\bar{\sP}(\vx_k, \mathbf{Y}_k)$ denotes the nominal joint distribution of $\vx_k$ and $\mathbf{Y}_k$ derived from the nominal system model~\eqref{eq1}. In the sense of the resulting minimum mean square error (MMSE), the optimal estimator can be found as  $\phi(\mathbf{Y}_k)=\E[\vx_k|\mathbf{Y}_k] \in \sH^\prime_{\mathbf{Y}_k}$. When $\bar{\sP}(\vx_k, \mathbf{Y}_k)$ follows a \textit{Gaussian distribution}, this solution corresponds to the standard ISE or Kalman filter, leveraging linearity and Gaussianity properties.

However, the problem in~\eqref{eq: DRE_optimization} becomes challenging when true joint state-measurement distribution $\sP(\vx_k, \mathbf{Y}_k)$ deviates from the nominal $\bar{\sP}(\vx_k, \mathbf{Y}_k)$, due to either uncertainties or outliers in the measurements. To address this challenge, the problem in~\eqref{eq: DRE_optimization} can be reformulated as
 \begin{equation} \label{eq: DRE_dual_min_max}
    \min_{\phi(\cdot) \in \sH^\prime_{{\mathbf{Y}}_{k}}} \max_{\sP(\vx_k, \mathbf{Y}_k) \in \sF} \E_{\sP(., .)} \big[\vx_k - \phi(\mathbf{Y}_k) \big]  \big[\vx_k - \phi(\mathbf{Y}_k) \big]^\top,
\end{equation} 
where $\sP(\vx_k, \mathbf{Y}_k) \in \sF$ denotes the true joint distribution of $\vx_k$ and $\mathbf{Y}_k$ constrained within the ambiguity set $\sF$, which defined as
\begin{equation*} 
 \sF = \left\{ \sP_{\xi} \left| \begin{aligned}
    &\xi \sim \sP_{\xi} \\
    &\sP_{\xi}(\xi \in \Xi) = 1 \\   
\end{aligned} \right. \right\}.
\end{equation*}
\begin{rem}
The ambiguity set is constructed as a neighborhood or "ball" around $\bar{\sP}_\xi$, allowing for some uncertainty while still assuming that the true distribution is reasonably close to the nominal one.
The size of this neighborhood, which is determined by a radius parameter, reflects the confidence level in the nominal distribution. Various approaches exist for defining this neighborhood, including metrics such as Kullback-Leibler divergence, $\tau$-divergence, $\phi$-divergence, and moment-based ambiguity sets.
Among these, the moment-based ambiguity sets are particularly attractive for their analytical tractability. It simplifies analysis by focusing on the first and second moments (mean and covariance), retaining robustness.
\end{rem}
The moment-based ambiguity set for the system states $\vx$ can be expressed as follows:  
\begin{equation*}
\sF_x = \left\{ \sP_x \left| \begin{aligned}
    & \vx \sim \sP_x \\  
    & \sP_x = \mathcal{N}(\vc_x, \mSigma^x) \\  
    & \sP_x(\vx \in \R^n) = 1 \\  
    & (\vc_x - \bar{\vx})^\top (\mP^{x})^{-1} (\vc_x - \bar{\vx}) \leq \theta_{3}^x \\  
    & \mSigma^x + (\vc_x - \bar{\vx})(\vc_x - \bar{\vx})^\top \leq \theta_{2}^x \mP^x \\  
    & \mSigma^x + (\vc_x - \bar{\vx})(\vc_x - \bar{\vx})^\top \geq \theta_{1}^x \mP^x  
\end{aligned} \right. \right\},
\end{equation*}
where $\bar{\vx}$ and $\mP^x$ represent the nominal mean and covariance of the state variable, and the parameters $\theta_1^x, \theta_2^x, \theta_3^x$ define the bounds of the set (\cite{DRE}). Similarly, the ambiguity set can be constructed for measurement noise $\vv$ as
\begin{equation*}
\sF_v = \left\{ \sP_v \left| \begin{aligned}
    & \vv \sim \sP_v \\  
    & \sP_v = \mathcal{N}(\vc_v, \mSigma^v) \\  
    & \sP_v(\vv \in \R^m) = 1 \\  
    & (\vc_v - 0)^\top \mR^{-1} (\vc_v - 0) \leq \theta_{3}^v \\  
    & \mSigma^v + (\vc_v - 0)(\vc_v - 0)^\top \leq \theta_{2}^v \mR \\  
    & \mSigma^v + (\vc_v - 0)(\vc_v - 0)^\top \geq \theta_{1}^v \mR  
\end{aligned} \right. \right\}.
\end{equation*} 
In this case, the nominal distribution of $\vv$ is assumed to be Gaussian with mean $\bm{0}$ and covariance $\mR$, i.e. $\overline{\sP}(\vv) = \mathcal{N}_m (\bm{0}, \mR)$. The parameters $\theta_1^v, \theta_2^v, \theta_3^v$ similarly define the bounds for the measurement noise ambiguity set. 
\begin{rem}
The choice of the moment-based ambiguity set is guided by practical considerations of uncertainty bounds. The assumed disturbance range aligns with common modeling practices, ensuring tractability while capturing distributional uncertainty. In real-world applications, suitable bounds can be identified through empirical data, allowing adaptation to specific scenarios.
\end{rem}
Since the focus is on applying state estimation to online recursive discrete-time systems, the problem is further reformulated into a one-step alternative
\begin{equation} \label{eq: DRE_dual_recursive}
    \min_{\phi(\cdot) \in \sH^\prime_y} \max_{\sP(\vx_k, \vy_k | {\vy}_{k-1}) \in \sF^\prime} \E_{\sP(., . | .)} \big[\vx_k - \phi(\vy_k) \big]  \big[\cdots]^\top,
\end{equation} 
where the new ambiguity set $\sF^\prime$ is constructed around the $\sP(\vx_k, \vy_k | \vy_{k-1})$ and the space of $\phi(\cdot)$ is only defined by previous measurement $\vy_k$ instead of $\mathbf{Y}_k$. 
This reformulation reduces the computational complexity of determining the optimal estimator at each time step. To solve~\eqref{eq: DRE_dual_recursive}, two key steps are necessary: \textit{i)} designing a suitable ambiguity set $\sF^\prime$ that accounts for parameter uncertainties and measurement outliers, and \textit{ii)} deriving explicit optimization equivalents for~\eqref{eq: DRE_dual_recursive} to enable efficient computation.  
Consequently, the study explores a distributionally robust Bayesian estimation problem 
 \begin{equation} \label{eq: DRE_primal}
    \min_{\phi(\cdot) \in \sH^\prime_{\vy}} \max_{\sP(\vx, \vy) \in \sF^{\prime\prime}} \E_{\sP(\vx_k, \vy_k)} \big[\vx_k - \phi(\vy_k) \big]  \big[\vx_k - \phi(\vy_k) \big]^\top,
\end{equation}  
subject to the nominal prior state distribution $\mathbb{\bar{P}}(\vx)$, the nominal conditional measurement distribution $\mathbb{\bar{P}}(\vy | \vx)$, and the true distribution $\sP(\vx, \vy)$ constrained in the ambiguity set $\sF^{\prime\prime}$. 
To solve the problem in~\eqref{eq: DRE_primal}, we are required to identify the least-favorable distribution from the ambiguity set $\sF^{\prime\prime}$. However, it depends on the specific choice of the estimator $\phi(\cdot)$. Therefore, we can alternatively try to solve the problem as
\begin{equation} \label{eq: DRE_dual_max_min}
\max_{\sP(\vx_k,\vy_k) \in \sF^{\prime\prime}} \min_{\phi(\cdot) \in \sH^{\prime}_\vy} \E_{\sP(\vx_k,\vy_k)}[\vx_k - \phi(\vy_k)][\vx_k - \phi(\vy_k)]^\top. 
\end{equation} 
 
To derive a robust estimation framework, we build upon the nominal distribution of $\vx\sim\mathcal{N}_n(\bar{\vx}, \mP^x)$ and $\vv\sim\mathcal{N}_m(\bm{0}, \mR)$, assuming $\vx$ and $\vv$ are independent. Let $\vs = \vy - \mC\bar{\vx}$ represent the innovation vector, with the covariance matrix of innovation $\mS$ defined as $\mS = \mC \mSigma^x \mC^\top + \mSigma^v$. The normalized innovation is denoted by $\vmu = \mS^{-1/2}\vs$. With the aforementioned definitions, the distributionally robust estimator can be formulated utilizing Theorem 9 in~\cite{DRE}. 

\section{Distributionally Robust Input and State Estimation} \label{sec: DRISE}

The ISE approach has been widely utilized to address the estimation challenges posed by the presence of unknown inputs and noisy measurements. In particular, ISE provides an optimal strategy for state estimation while simultaneously estimating the unknown inputs affecting the system.
However, classical ISE assumes Gaussian distributions for both process and measurement noises, limiting its robustness to deviations and outliers. To overcome this limitation, we propose the DRISE framework, which accounts for uncertainties in both distributions and covariances by employing moment-based ambiguity sets. The proposed approach effectively handles time-varying uncertainties, measurement outliers, and structural deviations in system dynamics.
The proposed DRISE can be summarized as follows.

Suppose the radius of the moment-based ambiguity sets are given by $\theta^{x} \geq 0$, $\theta_{2}^x \geq 1 \geq \theta_{1}^x \geq 0$, $\theta^{v} \geq 0$, $\theta_{2}^v \geq 1 \geq \theta_{1}^v \geq 0$, and at time $k$, with the nominal Gaussian prior conditional distribution of the state and the measurement noise 
\begin{equation*}
\bar{\sP}_{\vx_k|\vy_{k-1}} \sim \mathcal{N}_n(\hat{\vx}_{k}^\star, P_{k}^{x,\star}),\quad
\bar{\sP}_{\vv_k} = \mathcal{N}_m(\bm{0}, \mR_k),
\end{equation*}
then the distributionally robust state estimate $\hat{\vx}_{k}$ given $\vy_k$ is as follows.

\textit{Optimal Estimator:}
\begin{equation*} 
\hat{\vx}_{k} = \hat{\vx}_{k}^\star + \left(\mSigma_{k}^x \mC_k^\top - \mG_{k-1} \mM_k \mSigma^{v}_k \right) \mS_k^{-1/2} \cdot \psi \left(\mS_k^{-1/2} \vs_k \right),
\end{equation*}
where 
\begin{equation*}
\vs_k = \vy_k - \mC_k \hat{\vx}_{k}^\star, \quad \hat{\vx}_{k}^\star =  \hat{\vx}_{k}^- + \mG_{k-1} \hat{\vd}_{k-1},
\end{equation*}
and 
\begin{equation*}  
\mS_k = \mC_{k} \mSigma_{k}^x \mC_{k}^\top - \mC_{k} \mG_{k-1}\mM_{k} \mSigma_{k}^v - \mSigma_{k}^v \mM_{k}^\top \mG_{k-1}^\top \mC_{k}^\top + \mSigma_{k}^x.
\end{equation*}
The entry-wise function $\psi(\mu)$ for the $\epsilon$-contamination set, is defined as
\begin{equation*}  
\psi(\mu) = 
\begin{cases}
-K, & \mu \leq -K, \\
\mu, & |\mu| \leq K, \\
K, & \mu \geq K.
\end{cases}
\end{equation*}

\textit{Worst-Case Estimation Error Covariance:} 
\begin{equation*} 
\begin{aligned}
\hat{\mP}_{k}^x = \mSigma^{x}_k -  \left(\mSigma^{x}_k \mC_k^\top - \mG_{k-1} \mM_k \mSigma^{v}_k \right) \mS_k^{-1} \\
 \left(\mSigma^{x}_k \mC_k^\top - \mG_{k-1} \mM_k \mSigma^{v}_k \right)^{\top} \cdot i^{min}_{\mu},
\end{aligned}
\end{equation*}
where
 $i_{\mu}^{\min} = (1 - \epsilon)\left[1 - 2\Phi(-K)\right]$. The complete algorithm is presented in Algorithm~\ref{alg_DRISE}.

To demonstrate the validity of the proposed DRISE method, we establish its connection to the regular ISE approach presented in \cite{wan_attack_resilient}. By reformulating the optimal estimation equations and under specific assumptions, we show that DRISE reduces to ISE when the ambiguity sets are disregarded  ($\theta_2^x, \theta_2^v, i^{min}_{\mu} \to 1$) and the uncertainty models for the state and measurement noise distributions collapse to their nominal Gaussian distributions $(\mSigma_k^x \to \mP_k^x, \mSigma_k^v \to \mR_k)$. With considering $\psi \left(\mS_k^{-1/2} \vs_k \right)\to \mS_k^{-1/2} \vs_k $, where $\mS_k$ is the innovation matrix $\tilde{\mR}_k^\star$ in the ISE approach, defined as:
\begin{equation*}
\Scale[.94]{     
\tilde{\mR}_k^{\star} = \mC_k \mP_{k}^{x \star} \mC_k^\top  - \mC_k \mG_{k-1}\mM_k\mR_k-\mR_k\mM_k^\top \mG_{k-1}^\top \mC_k^\top+\mR_k.}
\end{equation*}
This reformulation highlights that DRISE extends ISE by incorporating distributionally robust optimization formulations to account for ambiguity sets in the state and measurement noise distributions.  

\begin{algorithm}[!t]  
\caption{Distributionally Robust Input and State Estimation: \\
$[\hat{\vd}_{k-1}, \mP^d_{k-1} \hat{\vx}_{k},\mP^x_{k}]= {\rm DRISE}(\hat{\vx}_{k-1},\mP^x_{k-1})$}
 \label{alg_DRISE} 
\begin{algorithmic}[1]    
 \FOR{$k = 1$ to $N$ }    
 
 		\STATEx \textbf{Step 1: Prediction:}  \label{Alg:Step1}
	    \STATE $\hat{\vx}_{k}^- = \mA_{k-1}\hat{\vx}_{k-1} + \mB_{k-1}\vu_{k-1}$;    		       
   		\STATE $\mP_{k}^{x,-} = \mA_{k-1}\mP^x_{k-1}\mA_{k-1}^\top + \mQ_{k-1}$;
        \STATE $\tilde{\mR}_k = \left(\mC_k\mP^{x,-}_{k}\mC_k^\top + \mR_k \right)^{-1}$;           
    
         \STATEx \textbf{  Step 2: Unknown input estimation:} \label{Alg:Step2}       
        \STATE $\mM_k = \left(\mG_{k-1}^\top \mC_k^\top \tilde{\mR_k}\mC_k \mG_{k-1} \right)^{-1} \mG_{k-1}^\top \mC_k^\top \tilde{\mR}_k$;      
        \STATE $ \hat{\vd}_{k-1} = \mM_k (\vy_k - \mC_k \hat{\vx}_{k}^-)$;        
    	\STATE $\mP^{d}_{k-1} = (\mG_{k-1}^\top \mC_k^\top \tilde{\mR}_k \mC_k \mG_{k-1})^{-1}$;
		\STATE $\mP^{x,d}_{k-1} = -\mP_{k-1}^x \mA_{k-1}^\top \mC_k^\top \mM_k^\top$;

  \STATEx  \textbf{  Step 3: Time Update:} \label{Alg:Step3}
        \STATE $\hat{\vx}_k^\star = \hat{\vx}_{k}^- + \mG_{k-1} \hat{\vd}_{k-1}$;        
         \STATE 
    $\begin{aligned}
               \mP^{x \star}_k =  \mA_{k-1}\mP^x_{k-1}\mA_{k-1}^\top + \mQ_{k-1} + \mG_{k-1} \mP^d_{k-1} \mG_{k-1}^\top \\
               + \mA_{k-1}\mP^{xd}_{k-1} \mG_{k-1}^\top + \mG_{k-1}(\mP_{k-1}^{xd})^\top \mA_{k-1}^\top  \\
               - \mG_{k-1}\mM_k \mC_k \mQ_{k-1} - \mQ_{k-1} \mC_k^\top  \mM_k^\top \mG_{k-1}^\top 
       \end{aligned}$

   \STATEx \textbf{  Step 4: Distributional Measurement update:} \label{Alg:Step4}
		 \STATE $\mSigma^{x}_k = \theta_{2}^x \mP^{x \star}_k \quad \text{and} \quad \mSigma^v_k = \theta_{2}^v \mR_k$;		
		 \STATE 
          $\begin{aligned}
          \mS_k =& \mC_k \mSigma^{x}_k \mC_k^\top - \mC_k \mG_{k-1}\mM_k \mSigma_{k}^v \\ -& \mSigma_{k}^v \mM_k^\top \mG_{k-1}^\top \mC_k^\top  + \mSigma_{k}^v;  \end{aligned}$
		 \STATE $\vs_{k} = \vy_k - \mC_k \hat{\vx}_k^\star$;			 
		 \STATE $\Scale[.94]{\hat{\vx}_{k} = \hat{\vx}_k^\star + \left(\mSigma^{x}_k \mC_k^\top - \mG_{k-1} \mM_k \mSigma^{v}_k \right)\mS_k^{-1/2} \psi_k (\mS_k^{-1/2} s_{k})}$;
		 \STATE 
		 $\begin{aligned}
		 \hat{\mP}^x_{k} = \mSigma^{x}_k - \left(\mSigma^{x}_k \mC_k^\top - \mG_{k-1} \mM_k \mSigma^{v}_k \right) \mS_k^{-1} \\
		  \left(\mSigma^{x}_k \mC_k^\top - \mG_{k-1} \mM_k \mSigma^{v}_k \right)^\top \cdot i^{min}_{\mu};
		 \end{aligned}$
\ENDFOR
   
\end{algorithmic}  
\end{algorithm}

\section{Experiments} \label{sec: exprements}

\begin{figure}
\begin{center}
\includegraphics[width=0.48\textwidth]{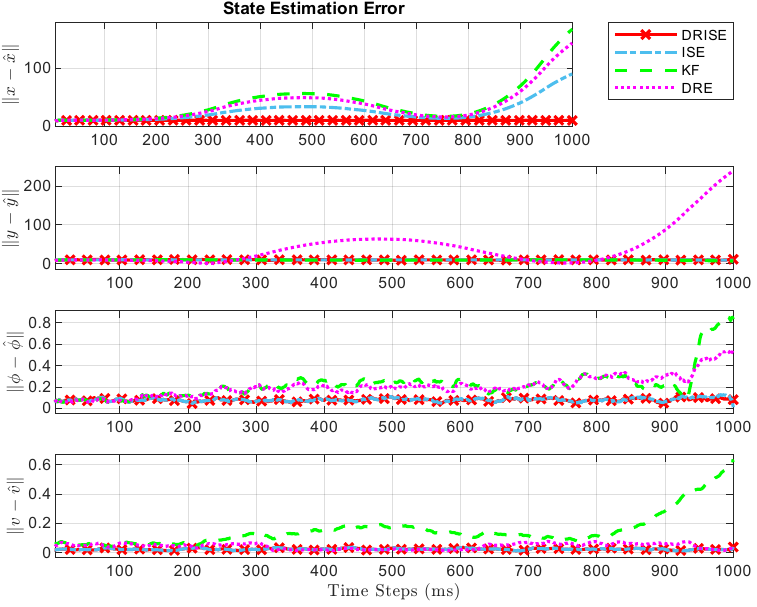} 
\caption{Comparison of state estimation errors.} 
\label{fig: tracking_error}
\end{center}
\end{figure}

We validate the proposed DRISE algorithm in the CARLA driving simulation environment (\cite{Dosovitskiy17}) and benchmark it against ISE, DRE, and the Kalman filter (KF) methods. The system parameters of the vehicle and the distributionally robust algorithm constants are selected based on~\cite{wan2021care} and ~\cite{DRE}, respectively.
The control input is defined as $\vu_k = [\beta_k, a_k]^\top$, where $\beta_k$ represents the vehicle’s slip angle and $a_k$ corresponds to the longitudinal acceleration. In particular, we have $\beta_k = \arctan\left(\frac{l_r}{l_f + l_r} \cdot \tan(\alpha_k) \right)$, where $\alpha_k$ represents the steering angle, $l_r$ and $l_f$ are the distances to the rear and front axles, respectively.
The unknown input vector is given as $\vd_k =  \left[\text{sign}\left(\sin(0.005 k)\right)\right] \left[1, 10 \right]^\top$,
where $k$ represents the discrete time step.
The noise covariances \(\mQ_k\) and \(\mR_k\) are represented as diagonal matrices where $\text{diag}(\mQ_k)=[1, 1, 0.01,0.001],$ and $\text{diag}(\mR_k)=[0.1, 0.1, 0.01, 0.001]$.
\begin{figure}[t]
\begin{center}
\includegraphics[width=0.48\textwidth]{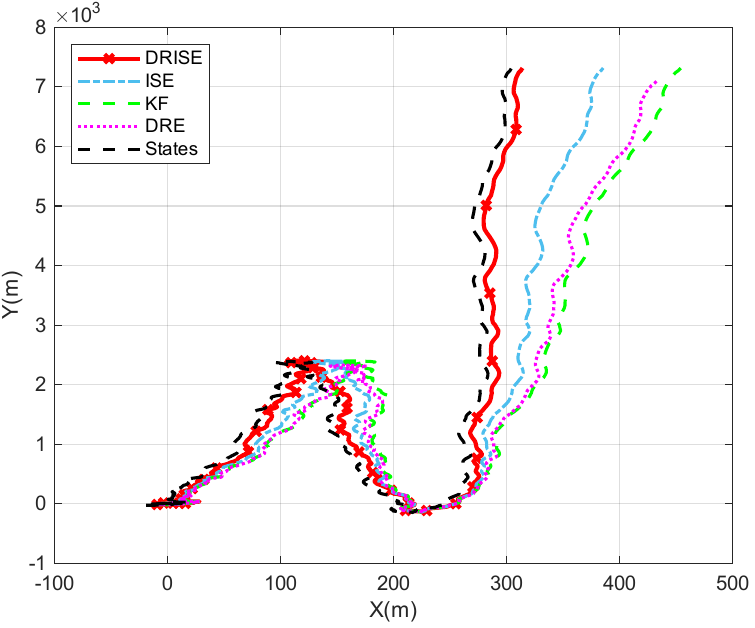}  
\caption{Trajectory tracking performance in X-Y plane.} 
\label{fig_2d}
\end{center}
\end{figure}

\begin{figure}
\begin{center}
\includegraphics[width=0.48\textwidth]{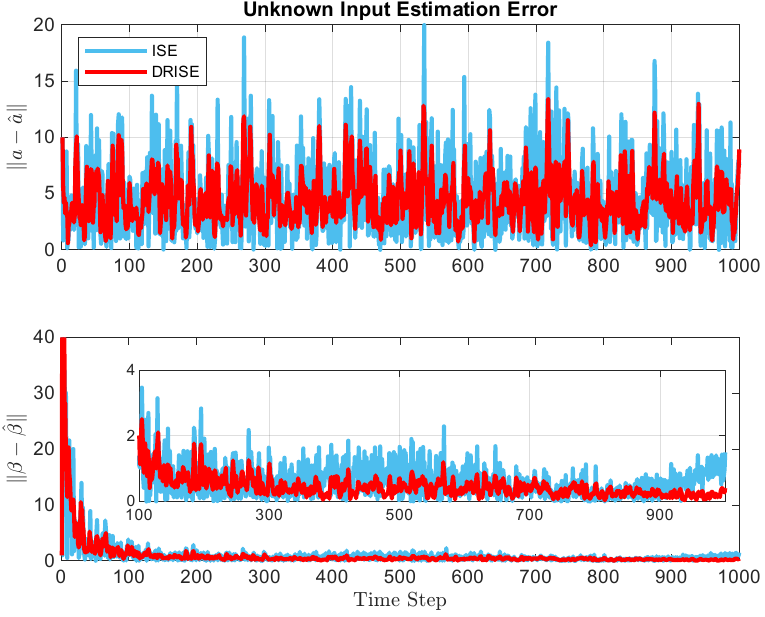}  
\caption{Comparison of unknown input estimation errors of DRISE and ISE.} 
\label{fig: attack_error}
\end{center}
\end{figure}

\begin{figure}
\begin{center}
\includegraphics[width=0.48\textwidth]{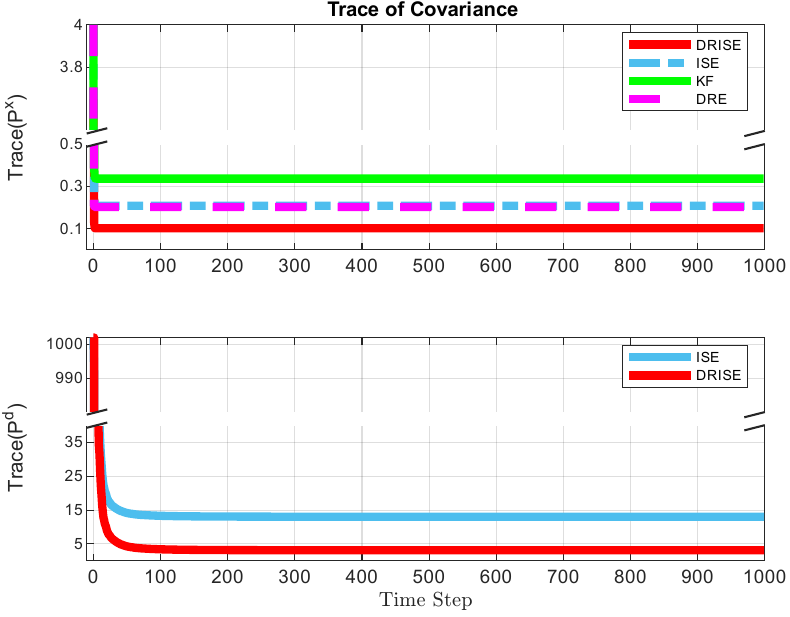}  
\caption{Trace of the unknown input estimation covariances.} 
\label{fig: Trace_px_pd}
\end{center}
\end{figure}

Figure~\ref{fig: tracking_error} shows the state estimation errors of all four algorithms. DRISE achieves the most accurate estimation. In contrast, ISE and DRE show moderate performance with slight drifts, while KF exhibits significant estimation errors, highlighting their limitations in handling uncertainties. 
The reference trajectory tracking comparison is illustrated in Figure~\ref{fig_2d}, where the vehicle is operated by the same baseline control law in CARLA under the same setting using four different estimation methods. The DRISE method follows the desired trajectory best, while others have a considerable drift.

The performance of unknown input estimation is compared in Figure~\ref{fig: attack_error}, which plots the estimation error of the acceleration $a_k$ and slip angle $\beta_k$ for DRISE and ISE. DRISE outperforms ISE, providing more accurate and reliable estimates of unknown inputs. This enhanced input estimation capability highlights the robustness of DRISE in scenarios where precise identification of external inputs is critical.
Figure~\ref{fig: Trace_px_pd} shows the trace of the covariance matrices for state and input estimation, and Table~\ref{fig: Table1} provides a quantitative comparison of the proposed DRISE algorithm with the ISE, DRE, and KF methods.

\begin{table}[!ht] 
\centering
\caption{Performance comparison.}
\begin{tabular}{c | c  c} 
 \hline
\textit{Method} & \textit{RMSE}($\hat{\vx}$) &  \textit{RMSE}($\hat{\vd}$) \\ \hline  \hline
\textbf{DRISE} & 14.21 & 7.48  \\ \hline
ISE   & 29.30 & 8.40 \\ \hline
DRE   & 47.37  &   -      \\ \hline
KF    & 69.23 &  -   \\ \hline
\end{tabular} \label{fig: Table1}
\end{table}

\section{Further Works} \label{sec: further}
This research sets the stage for incorporating more complex ambiguity sets as a foundational step in developing the DRISE framework. These alternative formulations will allow for a more versatile and tailored handling of uncertainty in various application scenarios. Additionally, future investigations will extend the evaluation of the DRISE framework to include complex dynamical systems to assess its robustness and performance under realistic, real-world conditions. This progression aims to validate the algorithm’s ability to effectively address uncertainties in diverse and challenging environments. 
Furthermore, future work will explore integrating alternative ambiguity sets, such as Wasserstein and KL-divergence-based sets, benchmarking their performance within the DRISE framework to assess their impact on robustness and estimation accuracy.

\section{Conclusion} \label{sec: conclusion}
This paper introduced the distributionally robust input and state estimation (DRISE) framework, which enabled simultaneous estimation of unknown inputs and states in uncertain environments.
By incorporating ambiguity sets and utilizing robust optimization techniques, DRISE provided robustness against model errors and measurement noises while maintaining computational efficiency. The recursive algorithm offered theoretical guarantees, including worst-case error bounds, making it suitable for adversarial scenarios and stochastic uncertainties, providing a reliable tool for complex and uncertain autonomous vehicles. 

\bibliography{ifacconf.bib}

\end{document}